\pdfoutput=1
\RequirePackage{ifpdf}
\ifpdf 
\documentclass[pdftex]{sigma}
\else
\documentclass{sigma}
\fi

\usepackage[mathscr]{eucal}

\newcommand{\mbf}[1]{\mathbf{#1}}
\newcommand{\msf}[1]{\text{\small $\sf{#1}$}}
\newcommand{\cmc}[1]{\mathcal{#1}}
\newcommand{\eus}[1]{\mathscr{#1}}
\newcommand{\euf}[1]{\mathfrak{#1}}
\newcommand{\bb}[1]{\mathbb{#1}}
\newcommand{\wh}[1]{\widehat{#1}}

\newcommand{\vt}{\vartheta}

\newcommand{\vp}{\varphi}
\newcommand{\om}{\omega}
\newcommand{\Om}{\Omega}

\newcommand{\AB}[1]{\langle#1\rangle}
\DeclareMathOperator{\id}{\normalfont\msf{id}}
\DeclareMathOperator{\ssid}{\text{\footnotesize\sf{id}}}

\newcommand{\C}{\bb{C}}

\newcommand{\E}{\bb{E}}

\newcommand{\N}{\bb{N}}

\newcommand{\R}{\bb{R}}

\newcommand{\cB}{\cmc{B}}

\newcommand{\sB}{\eus{B}}

\newcommand{\sS}{\eus{S}}

\newcommand{\eH}{\euf{H}}

\newcommand{\eS}{\euf{S}}

\newcommand{\U}{\mbf{1}}

\numberwithin{equation}{section}

\newtheorem{thm}{Theorem}[section]
\newtheorem{cor}[thm]{Corollary}
\newtheorem{lem}[thm]{Lemma}
\newtheorem{prop}[thm]{Proposition}
 { \theoremstyle{definition}
\newtheorem{emp}[thm]{}
\newtheorem{defi}[thm]{Definition}
\newtheorem{ex}[thm]{Example}
\newtheorem{ob}[thm]{Observation}
\newtheorem{rem}[thm]{Remark} }

\newcommand{\bemp}{\begin{emp}}
\newcommand{\eemp}{\end{emp}}
\newcommand{\bex}{\begin{ex}}
 \newcommand{\eex}{\end{ex}}
 \newcommand{\bexer}{\begin{exer}}
 \newcommand{\eexer}{\end{exer}}
 \newcommand{\bdefi}{\begin{defi}}
 \newcommand{\edefi}{\end{defi}}
 \newcommand{\brem}{\begin{rem}}
 \newcommand{\erem}{\end{rem}}
 \newcommand{\bob}{\begin{ob}}
 \newcommand{\eob}{\end{ob}}
 \newcommand{\bthm}{\begin{thm}}
 \newcommand{\ethm}{\end{thm}}
 \newcommand{\bprop}{\begin{prop}}
 \newcommand{\eprop}{\end{prop}}
 \newcommand{\bcor}{\begin{cor}}
 \newcommand{\ecor}{\end{cor}}
 \newcommand{\blem}{\begin{lem}}
 \newcommand{\elem}{\end{lem}}

\begin{document}
\allowdisplaybreaks

\newcommand{\arXivNumber}{0809.3538}

\renewcommand{\thefootnote}{}

\renewcommand{\PaperNumber}{071}

\FirstPageHeading

\ShortArticleName{Spatial Markov Semigroups Admit Hudson--Parthasarathy Dilations}

\ArticleName{Spatial Markov Semigroups Admit\\ Hudson--Parthasarathy Dilations\footnote{This paper is a~contribution to the Special Issue on Non-Commutative Algebra, Probability and Analysis in Action. The~full collection is available at \href{https://www.emis.de/journals/SIGMA/non-commutative-probability.html}{https://www.emis.de/journals/SIGMA/non-commutative-probability.html}}}

\Author{Michael SKEIDE}

\AuthorNameForHeading{M.~Skeide}

\Address{Universit\`a degli Studi del Molise, Dipartimento di Economia,\\Via de Sanctis, 86100 Campobasso, Italy}
\Email{\href{mailto:skeide@unimol.it}{skeide@unimol.it}}
\URLaddress{\url{ http://web.unimol.it/skeide/}}

\ArticleDates{Received February 23, 2022, in final form September 23, 2022; Published online October 03, 2022}

\Abstract{We present, for the first time, the result (from 2008) that (normal, strongly continuous) Markov semigroups on $\mathscr{B}(G)$ ($G$ a separable Hilbert space) admit a {\it Hudson--Parthasarathy dilation} (that is, a dilation to a cocycle perturbation of a {\it noise}) if and only if the Markov semigroup is {\it spatial} (that is, if it dominates an elementary CP-semigroup). The proof is by {\it general abstract nonsense} (taken from Arveson's classification of $E_0$-semigroups on $\mathscr{B}(H)$ by Arveson systems up to cocycle conjugacy) and not, as usual, by constructing the cocycle as a solution of a quantum stochastic differential equation in the sense of Hudson and Parthasarathy. All other results that make similar statements (especially, [\textit{Mem. Amer. Math. Soc.} \textbf{240} (2016), vi+126~pages, arXiv:0901.1798]) for more general $C^*$-algebras) have been proved later by suitable adaptations of the methods exposed here. (They use Hilbert module techniques, which we carefully avoid here in order to make the result available without any appeal to Hilbert modules.)}


\Classification{46L55; 46L53; 81S22; 60J25}

\renewcommand{\thefootnote}{\arabic{footnote}}
\setcounter{footnote}{0}

\section{Introduction}

(Quantum) {\it Markov semigroups} are models for irreversible evolutions of (quantum) physical systems. Dilating a Markov semigroup, means embedding the irreversible system into a reversible one in such a way that the original irreversible evolution can be recovered by projecting down (via a {\it conditional expectation}) the reversible evolution to the subsystem.

{\it Noises} are models for reversible systems containing a subsystem. A noise is actually a reversible evolution on a ``big'' system with a conditional expectation onto a ``small'' subsystem that leaves the small system invariant. One may think of a simultaneous description of a reversible system and the small system, but with the {\it interaction} switched off. When the interaction is switched on, the dynamics of the compound system changes and leaves the small system no longer invariant. The projection back to the small subsystem produces irreversible behavior.

Often, one tries to model the transition from the free dynamics (the noise) to the real dynamics by perturbation of the noise with a unitary cocycle. This is what we mean by a {\it Hudson--Parthasarathy dilation}. In practically all known examples, such cocycles have been obtained by means of a {\it quantum stochastic calculus}. (See Remark~\ref{HPhisrem}.) The {\it stochastic generator} of the cocycle, is composed from the generator of the Markov semigroup. Often, it may be interpreted in terms of an {\it interaction Hamiltonian}. In these notes, we show in the case of $\sB(G)$ (the algebra of all bounded operators on a Hilbert space~$G$) without using any calculus, that a Markov semigroup admits a Hudson--Parthasarathy dilation if and only if the Markov semigroup is {\it spatial}. The case of a Markov semigroup on a general von Neumann algebra (in particular, also on a commutative one, which corresponds to a classical dynamical system) is discussed, later, in Skeide~\cite{Ske16}.\looseness=-1

The necessary notions, spatial Markov semigroup, noise, Hudson--Parthasarathy dilation, and so forth, are explained in Section~\ref{notSEC} and then used to formulate the result. In Section \ref{proofSEC} we review the basic results about spatial $E_0$-semigroups and spatial product systems needed in the proof of the result, and we prove the result. In Section \ref{remSEC}, finally, we list several natural questions.

Our proof makes heavy use of Arveson's results \cite{Arv89} on the classification of $E_0$-semi\-groups (in particular, spatial ones) by tensor product systems of Hilbert spaces ({\it Arveson systems}). (We should emphasize that the main scope of \cite{Ske16} is not to generalize the present result to general von Neumann algebras -- which is easy --, but to fill a long standing gap, namely, to answer the question how Arveson's classification of $E_0$-semigroups by product systems generalizes to Hilbert modules.) We shall assume that the reader is familiar with Arveson's results at least in the spatial case, and we shall also assume that the reader knows the works by Bhat \cite{Bha96} and Arveson \cite{Arv97a} on the relation with Markov semigroups via the so-called {\it minimal weak dilation} (once more, in particular in the spatial case).

\section{Notations and statement of the result}\label{notSEC}

\bemp \label{EPS}\emph{$E_0$-semigroups}, in these notes,\footnote{By ``in these notes'' we refer to the fact that in most other of our papers we do not require any continuity of our semigroups. (Almost all constructions preserve strong continuity, in the sense that if we input a strongy continuous Markov (or CP) semigroup then the dilations we construct turn out to be strongly continuous, too. But continuity is not a necessary prerequisite for the constructions to work. And here, we require our Hilbert spaces separable, only because we wish to quote theorems from Arveson's theory.} are strongly continuous one-parameter semigroups of normal unital endomorphisms of $\sB(H)$ where $H$ is some separable infinite-dimensional Hilbert space.\footnote{Once for all other cases, \emph{strongly continuous} for a von Neumann algebra $\cB\subset\sB(H)$ refers to the strong operator topology (or the {\it point strong} topology). That is, a strongly continuous family $T_t$ of operators on $\cB$ satisfies that $t\mapsto T_t(b)h$ is norm continuous for all $b\in\cB$, $h\in H$.}
Arveson~\cite{Arv89} associated with every $E_0$-semigroup a tensor product system of Hilbert spaces (or, for short, an \emph{Arveson system}). He showed that two $E_0$-semigroups, $\vt^1$ on $\sB\big(H^1\big)$ and~$\vt^2$ on~$\sB\big(H^2\big)$ say, have isomorphic Arveson systems if and only if they are \emph{cocycle conjugate}. That is, there exist a unitary $u\colon H^1\rightarrow H^2$ and a strongly (and, therefore, $*$-strongly) continuous family of unitaries $u_t\in\sB\big(H^2\big)$ fulfilling:
\begin{enumerate}\itemsep=0pt
\item
The $u_t$ form a \emph{left cocycle} with respect to $\vt^2$, that is, $u_{s+t}=u_s\vt^2_s(u_t)$ for all $s,t\in\R_+$.

\item $u\vt^1_t(u^*au)u^*=u_t\vt^2_t(a)u_t^*$ for all $t\in\R_+$, $a\in\sB\big(H^2\big)$.
\end{enumerate}
\eemp

\bemp \label{WD}
\emph{Markov semigroups}, in these notes, are strongly continuous one-parameter semigroups of normal unital completely positive maps (\emph{CP-maps}) on $\sB(G)$ where $G$ is some separable Hilbert space. Bhat \cite[Theorem~4.7]{Bha96} states that every Markov semigroup~$T$ on $\sB(G)$ admits a \emph{weak dilation} to an $E_0$-semigroup $\vt$ on $\sB(H)$ in the following sense: There is an isometry $\xi\colon G\rightarrow H$ such that
\[
T_t(b)=\xi^*\vt_t(\xi b\xi^*)\xi
\]
for all $t\in\R_+,b\in\sB(G)$.\footnote{Frequently, $G$ is identified with the subspace $\xi G$ of $H$, so that $\sB(G)$ would be identified with the corner $\xi\,\sB(G)\xi^*=p\sB(H)p$ of $\sB(H)$, where $p$ denotes the projection in $\sB(H)$ onto $\xi G$. That would lead to slightly more readable formulae. But, in these notes, we will, have -- by definition (see Section~\ref{noise}) -- sitting $\sB(G)$ as a unital subalgebra of $\sB(H)$ so that $\id_H\in\sB(G)\subset\sB(H)$ needs to be carefully distinguished from the corner $p\sB(H)p$ (of course, isomorphic to~$\sB(G)$).}
This weak dilation can be chosen \emph{minimal} in the sense that $H$ has no proper subspace containing $\xi G$ and invariant under $\vt_{\R_+}(\xi\,\sB(G)\xi^*)$. The minimal weak dilation is unique up to suitable unitary equivalence, namely, by a unique unitary that sends one $\xi$ to the other. Bhat \cite{Bha96} defines the Arveson system \emph{associated} with the Markov semigroup to be the Arveson system of its minimal weak dilation $\vt$.\footnote{Strictly speaking, for an $E_0$-semigroup as defined in the first paragraph, $H$ should be infinite-dimensional. However, the missing case, where the $H$ of the minimal weak dilation is finite-dimensional, happens if and only if $G$ is finite-dimensional and $T$ consists of automorphisms. In this case, for our theorem below there is nothing to prove, and we will tacitly exclude it from the subsequent discussion.} Here we need only know that the Arveson system of any weak dilation of $T$ contains the Arveson system of $T$.
\eemp

\bemp \label{noise} By a \emph{noise over $\sB(G)$} we understand an $E_0$-semigroup $\sS$ on some $\sB(H)$ that contains $\sB(G)$ as a unital $W^*$-subalgebra\footnote{That is, we identify elements of $\sB(G)$ with their image under a(n unnamed) normal unital representation $\sB(G)\rightarrow\sB(H)$.}
and an isometry $\om\colon G\rightarrow H$, fulfilling the following:
\begin{enumerate}\itemsep=0pt
\item $\sS$ leaves $\sB(G)$ \emph{invariant}, that is, $\sS_t(b)=b$ for all $t\in\R_+,b\in\sB(G)\subset\sB(H)$.

\item $\om$ is $\sB(G)$-$\C$-bilinear, that is, $b\om=\om b$ for all $b\in\sB(G)$.

\item $\sS$ leaves $\E\colon a\mapsto\om^*a\om$ \emph{invariant}, that is, $\E\circ\sS_t=\E$ for all $t\in\R_+$.
\end{enumerate}
By (2), $\E$ is a conditional expectation (onto $\sB(G)$ sitting as unital subalgebra of $\sB(H)$). By~(3) (applied to $a=\om b\om^*$), $\sS$ with the isometry $\om$ is a weak dilation of the identity Markov semigroup on $\sB(G)$. Therefore, the projection $p:=\om\om^*\in\sB(H)$ is \emph{increasing}, that is, $\sS_t(p)\ge p$ for all $t\in\R_+$. A noise is \emph{reversible}, if all $\sS_t$ are automorphisms. In this case, $\sS_t(p)=p$ for all $t\in\R_+$.\footnote{In fact, suppose $\sS$ is implemented as $\sS_t=v_t\bullet v_t^*$ by a unitary semigroup $v_t$. One easily checks that $p$ is increasing if and only if $\om^*v_t^*\om$ is an isometry. Since $\sS_t(b)=b$ we have $bv_t=v_tb$. Since also $b\om=\om b$, it follows that $\om^*v_t^*\om$ is in the center of $\sB(G)$ and, therefore, a unitary. So, $p$ is also decreasing, thus, constant.\label{revcFN}}
\eemp

\brem This definition of noise is from Skeide \cite{Ske06d}. In the scalar case (that is, $G=\C$) it corresponds to noises in the sense of Tsirelson \cite{Tsi98p1,Tsi03p1}. A reversible noise is close to a {\it Bernoulli shift} in the sense of Hellmich, K\"ostler and K\"ummerer \cite{HKK04p}.
\erem

\bdefi
A \emph{Hudson--Parthasarathy dilation} (\emph{HP-dilation}, for short) of a Markov semigroup $T$ on $\sB(G)$ is a noise $(\sS,\om)$ and a unitary left cocycle $u_t$ with respect to $\sS$, such that the cocycle conjugate $E_0$-semi\-group $\vt$ defined by setting $\vt_t(a)=u_t\sS_t(a)u_t^*$, fulfills
\[
\E\circ\vt_t\upharpoonright\sB(G)=T_t
\]
for all $t\in\R_+$.

The HP-dilation is \emph{reversible}, if the underlying noise is reversible. In this case, also $\vt$ is an automorphism semigroup.

The HP-dilation is \emph{weak} if $\vt$ with the isometry $\om$ is also a weak dilation of $T$.
\edefi

\bob \label{wob}
Recall that being a weak HP-dilation is stronger a condition than being an HP-dilation. Since reversible noises have the trivial (that is, the one-dimensional) Arveson system, not much remains for what a Markov semigroup $T$ can be if it admits a weak and reversible HP-dilation. (One can show, that a Markov semigroup with trivial Arveson system has the form $T_t=w_t^*\bullet v_t$ for a semigroup $v_t$ of isometries in $\sB(G)$.\footnote{This is a statement that follows most easily using Hilbert (or, better, von Neumann) modules. (The GNS-product system of von Neumann $\sB(G)$-correspondences of a normal CP-semigroup on $\sB(G)$ is $\sB(G)\otimes\eH^\otimes$ where~$\eH^\otimes$ is the Arveson system of that CP-semigroup; this first in Bhat and Skeide \cite[Section 13]{BhSk00} though there not very explicit. If $\eH^\otimes$ is the trivial Arveson system, then the units of the GNS-product system are just the semigroups in~$\sB(G)$. Conversely, a CP-semigroup of the form $c_t^*\bullet c_t$ has the trivial product system of $\sB(G)$-correspondences as GNS-product system (see Shalit and Skeide \cite[Observation~7.2]{ShaSk10p}) and, therefore, the trivial Arveson system.) The argument does not fit into this note that avoids module techniques, so we omit a formal proof.}
Taking any unitary dilation $u_t\in\sB(H)$ and the trivial (=~constant) noise on $\sB(H)$ (so that the unitary semigroup is also a cocycle) we get an HP-dilation.) However, since any $E_0$-semigroup on $\sB(H)$ extends to an automorphism semigroup on some $\sB\big(\wh{H}\big)\supset\sB(H)\ni\id_{\wh{H}}$ (Arveson and Kishimoto~\cite{ArKi92} or Skeide~\cite{Ske07a}), existence of a weak HP-dilation grants existence of a reversible HP-dilation -- a reversible HP-dilation having {\it sitting inside} a weak HP-dilation in a specific way. In general, it is still unknown whether existence of a reversible or general HP-dilation implies existence of a weak HP-dilation. In particular, since by our main result, Theorem~\ref{mthm} below, spatial Markov semigroups admit weak HP-dilations, it is unknown whether there are non-spatial Markov semigroups that admit other types of HP-dilations. So, the {\it only-if} part of the theorem might fail, if we drop {\it weak} from the hypotheses on the HP-dilation.
\eob

\brem \label{HPhisrem}
Since the seminal work of Hudson and Parthasarathy \cite{HuPa84,HuPa84a}, the cocycles of Hudson--Parthasarathy dilations have been obtained with the help of {\it quantum stochastic calculi} as solutions of {\it quantum stochastic differential equations}. \cite{HuPa84} dealt with a {\it Lindblad generator} with {\it finite degree of freedom}, while \cite{HuPa84a} considers a general (bounded) Lindblad generator. Chebotarev and Fagnola \cite{CheFa98} deal with a large class of unbounded generators. Versions for general von Neumann algebras (Goswami and Sinha~\cite{GoSi99}, K\"ostler \cite{Koes00}) or $C^*$-algebras (Skeide \cite{Ske00}) require Hilbert modules. (Apart from the fundamental monograph \cite{Par92} by Parthasarathy, a still up-to-date reference for everything that has to do with calculus based on Boson Fock spaces or modules are Lindsay's lecture notes~\cite{Lin05}. Results that use other types of Fock constructions or abstract representation spaces are scattered over the literature.)

We get our cocycle in an entirely different way, appealing to granted cocycle conjugacy of $E_0$-semigroups having the same Arveson system as discussed in Section~\ref{EPS}.
\erem

To, finally, formulate our main result, we need a last notion.

\bdefi[{Arveson \cite[Definition 2.1]{Arv97a}}]
A \emph{unit} for a Markov semigroup $T$ on $\sB(G)$ is a~strongly continuous semigroup $c$ in $\sB(G)$ such that $T$ \emph{dominates} the elementary CP-semigroup $S_t(b):=c_t^*bc_t$, that is, the difference $T_t-S_t$ is completely positive for all $t\in\R_+$. A Markov semigroup on $\sB(G)$ is \emph{spatial}, if it admits units.
\edefi

And, now, here is our main result:

\bthm \label{mthm}
A Markov semigroup on $\sB(G)$ is spatial if and only if it admits a weak Hudson--Partha\-sarathy dilation. Such a Hudson--Parthasarathy dilation may be extended to a reversible Hudson--Parthasarathy dilation.
\ethm

\section{Proof}\label{proofSEC}

We said, a Markov semigroup is spatial if it admits a unit. Of course, also an $E_0$-semigroup is a Markov semigroup. For $E_0$-semigroups there is Powers' definition \cite{Pow87} of spatiality in terms of {\it intertwining semigroups of isometries}, also referred to as (isometric) units. Also for Arveson systems there is the concept of {\it units} and an Arveson system is \emph{spatial}, if it admits a unit. We do not repeat here the definition of Arveson system nor that of a unit for an Arveson system, but simply will put together known statements about them. (The discussion for general von Neumann algebras in Skeide~\cite{Ske16} is much more self-contained, and many of the statements for~$\sB(G)$, we simply quote here, will drop out there very naturally without any effort.) Bhat~\cite[Section~6]{Bha01} compared several notions of units. What is important to us is that all concepts of spatiality in the sense of existence of units coincide: A semigroup, Markov or~$E_0$, is spatial in whatsoever sense if and only if its associated Arveson system is spatial; moreover, whether a weak dilation is spatial or not, does not depend on whether the dilation is minimal, but only on whether the dilated Markov semigroup is spatial or not; see Bhat \cite[Section~6]{Bha01} or Arveson \cite[Sections~8.9 and~8.10]{Arv03}.

Let us start with the more obvious parts of the theorem. It is the word {\it weak} that guarantees that the Markov semigroup that admits a weak HP-dilation has to be spatial. (A~noise has a~spatial Arveson system, so that also the Arveson system of the cocycle conjugate HP-dilation $\vt$ is spatial. Since $\vt$ is also a weak dilation, by the discussion in Section~\ref{WD}, its Arveson system contains the Arveson system of~$T$. This is enough to see that $T$ is spatial; see Bhat \cite[Section~6]{Bha01} or, more recently, Bhat, Liebscher, and Skeide~\cite{BLS10}, but see also the discussion below of the relation between noises over $\sB(G)$ and noises over~$\C$.) We mentioned already in Observation~\ref{wob} that a weak HP-dilation {\it embeds} into a reversible one.

So, what remains to be shown is the construction of a weak HP-dilation for any spatial Markov semigroup. Therefore, for the balance, let us fix a spatial Markov semigroup $T$.

Every weak dilation of $T$ is spatial. This has important consequences: For a spatial Arveson system it is easy to construct an $E_0$-semigroup $\eS$ on some $\sB(\eH)$ that has as associated Arveson system the one we started with; see \cite[Appendix]{Arv89}. Moreover, it is easy to see that for this $E_0$-semigroup there exists a unit vector $\Om\in\eH$ such that $\eS$ leaves the state $\vp:=\AB{\Om,\bullet\Om}$ \emph{invariant} (that is, $\vp\circ\eS_t=\vp$ for all $t\in\R_+$); in other words, $\eS$ is a noise over $\C$.%
\footnote{Arveson's construction is the core of any forthcoming construction of an $E_0$-semigroup from a unital (and {\it central} for modules) unit in a product system, by an inductive limit. It is easy to see that the elements of the unit are, under the inductive limit, sent to a fixed unit vector $\Om$ and that this $\Om$ does the job. Moreover, vector states are pure, so the $\eS$ is an $E_0$-semigroup in {\it standard form}; see Powers \cite{Pow99} and Alevras \cite{Ale01}, but also Bhat and Skeide \cite{BhSk00,Ske02,Ske06d} for the module case.}
Two $E_0$-semigroups are cocycle conjugate if and only they have isomorphic Arveson systems; see the corollary of Arveson \cite[Definition 3.20]{Arv89}. So, a Markov semigroup is spatial if and only if one (and, therefore, all) weak dilation(s) is (are) cocycle conjugate to an $E_0$-semigroup with an invariant vector state.

From $E_0$-semigroups with invariant vector states to noises and back, there is only a small step. Suppose we have a noise $(\sS,\om)$ over $\sB(G)$ on $\sB(H)$. Clearly, a unital $W^*$-subalgebra $\sB(G)$ of $\sB(H)$ decomposes $H$ into $H$``$=$''$G\otimes\eH$ for some multiplicity space $\eH$, and $\sS$ leaves $\sB(G)$``$=$''$\sB(G)\otimes\id_\eH$ invariant if and only if $\sS_t=\id_{\sB(G)}\otimes\eS_t$ for a unique $E_0$-semigroup $\eS$ on $\sB(\eH)$.\footnote{$\sS_t(\ssid_G\otimes a)$ is in the commutant of $\sS_t(\sB(G))=\sB(G)$. So $\sS_t$ leaves $\ssid_G\otimes\sB(\eH)$ invariant.}
For that the isometry $\om$ intertwines the actions of $\sB(G)$, it necessarily has the form $\om=\id_G\otimes\Om\colon g\mapsto g\otimes\Om$ for a unique unit vector $\Om\in\eH$. Clearly, $\sS$ leaves the conditional expectation $\E$ invariant if and only if $\eS$ leaves the vector state $\vp:=\AB{\Om,\bullet\Om}$ invariant. Conversely, if $\eS$ is an $E_0$-semigroup with an invariant vector state $\vp$ induced by a unit vector $\Om\in\eH$, then the $E_0$-semigroup $\sS=\id_{\sB(G)}\otimes\,\eS$ on $\sB(G\otimes\eH)$ with the isometry $\om:=\id_G\otimes \,\Om$ is a noise over~$\sB(G)$. On the other hand, $\sS$ is just a multiple of $\eS$, and multiplicity does not change the Arveson system; see Arveson \cite[Poposition~3.15]{Arv89}. Therefore, the Arveson system of a noise $\sS$ is spatial. Consequently, a Markov semigroup is spatial if and only if one (and, therefore, all) weak dilation(s) is (are) cocycle conjugate to a noise.

To construct a weak Hudson--Parthasarathy dilation for $T$, we need to choose the noise (to which a weak dilation is cocycle conjugate) in such a way that it admits a Hudson--Parthasarathy cocycle (which turns the noise by cocycle perturbation into a dilation). The Arveson system of the weak dilation $\vt$ on $\sB(H)$ of $T$ with the isometry $\xi$ is spatial. To that Arveson system construct an $E_0$-semigroup $\eS$ on a Hilbert space $\eH$ with an invariant vector state $\vp=\AB{\Om,\bullet\Om}$. We tensor it with the identity on $\sB(G)$ as described before to obtain a noise $(\sS,\om)$ with the same Arveson system as $\vt$. We wish to identify the two Hilbert spaces (infinite-dimensional and separable, unless $T$ is an automorphism semigroup on $M_n$) by a unitary $u\colon H\rightarrow G\otimes\eH$ in such a way that $u\xi=\om$. But this is easy. If $T$ is an $E_0$-semigroup, then, since $T$ is its own weak dilation, there is nothing to show. If~$T$ is not an $E_0$-semigroup, then both $E_0$-semigroups, $\vt$~and~$\eS$, are proper. We simply fix a unitary $u\colon H\rightarrow G\otimes\eH$ that takes $\xi g$ to $g\otimes\Om=\om g$ and is arbitrary on the (infinite-dimensional!\footnote{Clear is that the scalar noise $\eS$ lives on an infinite-dimensional $\eH$ so that $\infty=\dim\Om^\perp\le\dim(G\otimes\Om)^\perp$. It is a bit of {\it folklore} that also $(\xi G)^\perp$ in a weak dilation of a non-$E_0$-Markov semigroup has to be infinite-dimensional, but it is not so easy to find a reference for that. It is a bit inconvenient that for seeing that, we would have to enter the construction of the minimal dilation (contained in any other weak dilation). It is enough to mention that the construction of the minimal dilation as in Bhat and Skeide \cite[Section 5]{BhSk00} (for modules but only over~$\sB(G)$, for which the ``translation'' in \cite[Section 13]{BhSk00} can be applied) is also (unlike Bhat's in \cite{Bha96}!) by an inductive limit, and that it is easy to see that in each step of the inductive limit ``new space'' is added, if the Markov semigroup is nonendomorphic. (By `step', we mean that we restrict the parameter set $\R_+$ to the discrete case $\N_0\subset\R_+$ and the isometry embedding the space at time $n$ into that at time $n+1$ is always proper. (The latter part is a~retraction of the fact that the Stinespgring isometry of the Stingespring construction of a~proper Markov map is always a~proper isometry. This implies the former part, because the isometry going from $n$ to $n+1$ is (for $\sB(G)$, at least!) just an amplification of the Stinestpring isometry.)}) complements of $\xi G$ and $G\otimes\Om$. There exists, then, a~left cocycle $u_t$ with respect to $\sS$ that fulfills
\[
u\vt_t(u^*au)u^*=u_t\sS_t(a)u_t^*.
\]
We find
\begin{align*}
T_t(b) & = \xi^*\vt_t(\xi b\xi^*)\xi=\xi^*u^*u\vt_t(u^*u\xi b\xi^*u^*u)u^*u\xi\\
& = \om^*u\vt_t(u^*\om b\om^*u)u^*\om=\om^*u_t\sS_t(\om b\om^*)u_t^*\om,
\end{align*}
so that $u_t\sS_t(\bullet)u_t^*$ with the isometry $\om$ is a weak dilation of~$T$. In particular, the projection $\om\om^*$ must be increasing, that is, $u_t\sS_t(\om\om^*)u_t^*\om\om^*=\om\om^*$ or $u_t\sS_t(\om\om^*)u_t^*\om=\om$. Now, by the special property of $\om$, we have $\om b\om^*=(\om\om^*)b(\om\om^*)$. It follows
\[
T_t(b)
=\om^*u_t\sS_t(\om \om^*)u_t^*u_t\sS_t(b)u_t^*u_t\sS_t(\om \om^*)u_t^*\om
=\om^*u_t\sS_t(b)u_t^*\om,
\]
that is, the cocycle perturbation of the noise $(\sS,\om)$ by the cocycle $u_t$ is a Hudson--Parthasarathy dilation of $T$.

\section{Remarks and outlook}\label{remSEC}

Our construction of a Hudson--Parthasarathy dilation is by completely abstract means. This leaves us with a bunch of natural questions.

Is our cocycle in any way adapted? Hudson--Parthasarathy cocycles obtained with quantum stochastic calculus on the Boson Fock space are adapted in the sense that~$u_t$ is in the commutant of $\sS_t(\U\otimes\sB(\eH))$ for each $t\in\R_+$. Is our cocycle possibly adapted in the sense of \cite[Definition~7.4]{BhSk00}? (Today, we would prefer to say weakly adapted. Roughly, this means $\sS_t(\om\om^*)u_t^*\om\om^*=u_t^*\om\om^*$, so that the $u_t^*\om\om^*$ form a partially isometric cocycle.) Instead of the minimal Arveson system of the minimal weak dilation, we could have started the constructive part with the Arveson system associated with the free flow {\it generated} by the spatial minimal Arveson system in a sense to be worked out in Skeide~\cite{Ske08p2}. (This has been outlined in Skeide~\cite{Ske06d}.) These free flows come along with their own notion of adaptedness (see K\"ummerer and Speicher \cite{KueSp92}, Fowler~\cite{Fow95}, and Skeide~\cite{Ske00}), and we may ask whether the cocycle is adapted in this sense.

In any case, quantum stochastic calculi, also the abstract one in K\"ostler~\cite{Koes00}, provide a relation between additive cocycles and multiplicative (unitary) cocycles. Differentials of additive cocycles are, roughly, the differentials of the quantum stochastic differential equation to be resolved. We may ask, whether this relation holds for all spatial Markov semigroups, also if they are not realized on the Fock spaces, that is, on {\it type I} or {\it completely spatial} noises. The additive cocycles, usually, take their ingredients from the generator of the Markov semigroup. If that generator is bounded, then one may recognize the constituents of the {\it Christensen--Evans} generator, or, in the $\sB(G)$-case, of the {\it Lindblad} generator. This raises the problem to characterize spatial Markov semigroups in terms of their generators. Do they have generators that resemble in some sense the Lindblad form? Apparently the most general form of unbounded Lindblad type generators of Markov semigroups on $\sB(G)$ has been discussed in Chebotarev and Fagnola~\cite{CheFa98}. But only a subclass of these generators could be dilated by using Hudson--Parthasarathy calculus on the Boson Fock space. It seems natural to expect that Markov semigroups having this sort of generators are all spatial. Is it possibly that they could not be dilated on the Fock space because their product systems are type II (non-Fock) and not completely spatial (Fock)? In this case, can the solution be obtained with a calculus on noises emerging from type II product systems that are no longer CCR-flows (Fock noises)? In any case, whenever for an example a solution of the problem has been obtained with calculus, then we may ask, whether our abstract cocycle (which, of course, can be written down explicitly; see Skeide~\cite{Ske16}) can be related to the concrete cocycle emerging from calculus. Generally, we may ask, how two possible cocycles with respect to the same noise (adapted in some sense or not) are related. (This question connects to properties of the automorphism group of spatial Arveson systems such as discussed in Tsirelson~\cite{Tsi08}.)

Last but not least, we ask, if there exist nontrivial examples of nonspatial Markov semigroups. By this we mean Markov semigroups that are not type~III $E_0$-semigroups, or tensor products of such with a spatial Markov semigroup. A big step towards answering this question has been taken in Skeide~\cite{Ske12}, where (starting from our construction Skeide~\cite{Ske06} of an $E_0$-semigroup for every Arveson system) we constructed for each Arveson system a non-endomorphic Markov semigroup. If we find a type III Arveson system that does not factor into a tensor product with a spatial product system, then we would find a properly non-spatial Markov semigroup.

\subsection*{Acknowledgements}

I would like to thank Rajarama Bhat, Malte Gerhold, Martin Lindsay, and Orr Shalit for several valuable comments and bibliographical hints, which improved the final version in several ways. In particular, I wish to express my gratitude to the referees. One of them spotted a seriously unclear point in a definition and a small gap in an argument. This work was supported by research funds of University of Molise and Italian MIUR under PRIN 2007.

\pdfbookmark[1]{References}{ref}
\LastPageEnding

\end{document}